\documentclass{proc-l}

\usepackage{amssymb}
\usepackage{graphicx}
\usepackage{comment}
\usepackage[shortlabels]{enumitem}

\newcounter{theo}
\newcounter{qstn}
\newcounter{dfn}
\newcounter{rmk}
\newtheorem{defn}[dfn]{Definition}
\newtheorem{lem}[theo]{Lemma}
\newtheorem{thm}[theo]{Theorem}
\newtheorem{qstin}[qstn]{Question}
\newtheorem{cor}[theo]{Corollary}
\newtheorem{prop}[theo]{Proposition}
\newtheorem{fact}[theo]{Fact}

\newtheoremstyle{nopunct}
{\topsep}{\topsep}
{\itshape}{}
{\bfseries}{}{.5em}{}
\theoremstyle{nopunct}
\newtheorem*{thm*}{Theorem}
\newtheorem*{prop*}{Proposition}
\newtheorem*{cor*}{Corollary}

\theoremstyle{remark}
\newtheorem{remark}[rmk]{Remark}

\newcommand{\sJ}{\mathcal{J}}
\newcommand{\sF}{\mathcal{F}}
\newcommand{\sR}{\mathcal{R}}
\newcommand{\sB}{\mathcal{B}}
\newcommand{\sN}{\mathcal{N}}
\newcommand{\sM}{\mathcal{M}}
\newcommand{\sG}{\mathcal{G}}

\newcommand{\A}{\mathbb{A}}
\newcommand{\R}{\mathbb{R}}
\newcommand{\Pn}{\mathbb{P}}
\newcommand{\Q}{\mathbb{Q}}

\newcommand{\model}{\vDash}
\newcommand{\force}{\Vdash}

\newcommand{\cf}{\text{cf}}
\newcommand{\cc}{\mathfrak{c}}
\newcommand{\bcl}{\backslash}
\newcommand{\ol}[1]{\overline{#1}}
\newcommand{\id}{\textsf{id}}
\newcommand{\ran}{\text{ran }}
\newcommand{\dom}{\text{dom }}

\newcommand{\thry}[1]{\textsf{#1}}

\begin{document}

\title{Independent families under determinacy-like hypotheses}

\author[S. Bahrke]{Stewart Bahrke$^\dagger$}
\address{Department of Mathematics, University of North Texas, Denton, TX, USA}
\email{bahrke.stewart@gmail.com}
\thanks{
$^\dagger$This work was partially supported by the Simons Foundation grant (award no. SFI-MPS-T-Institutes-00010825) and from State Treasury funds as part of a task commissioned by the Minister of Science and Higher Education under the project “Organization of the Simons Semesters at the Banach Center - New Energies in 2026-2028” (agreement no. MNiSW/2025/DAP/491).
\includegraphics[scale=.33]{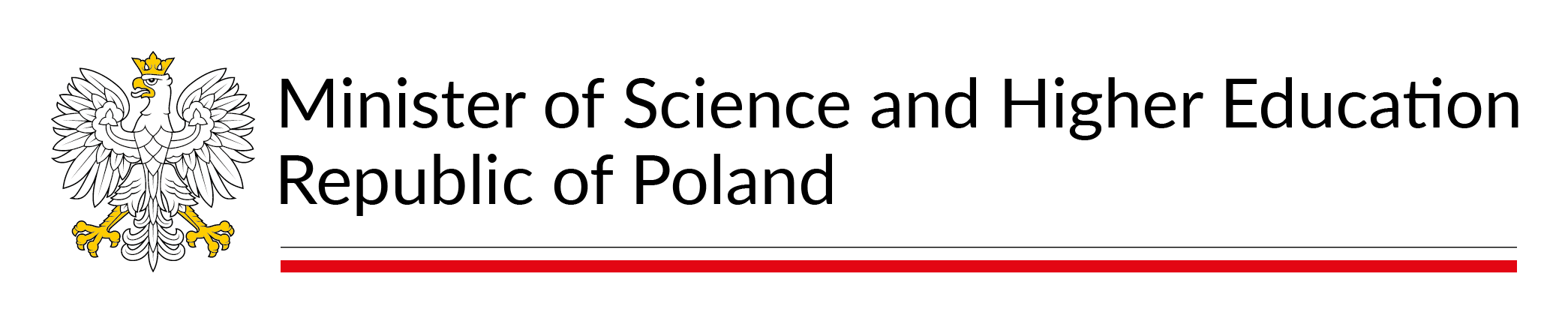}
\includegraphics[scale=.33]{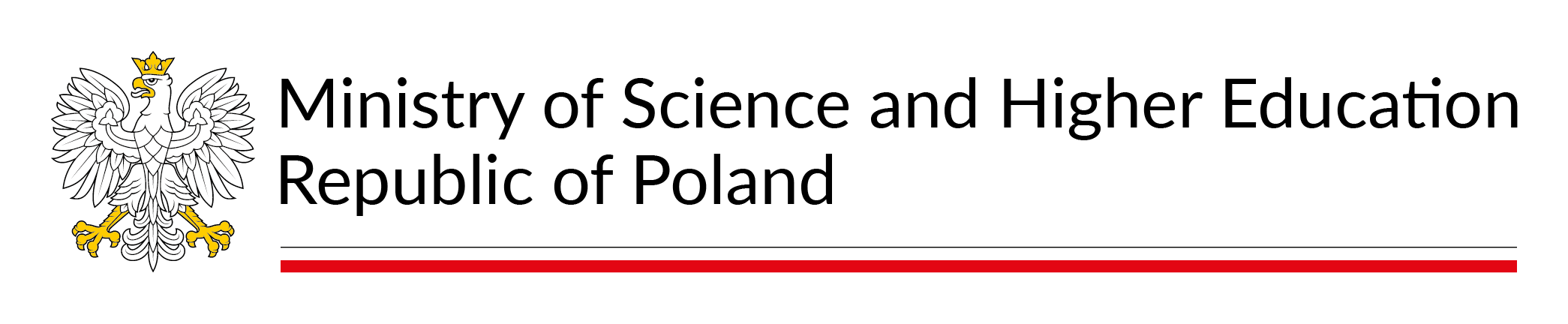}}

\author[N. Trang]{Nam Trang$^\ddagger$}
\address{Department of Mathematics, University of North Texas, Denton, TX, USA}
\email{Nam.Trang@UNT.edu}
\thanks{$^\ddagger$The second author acknowledges support by the National Science Foundation under Career Award No. 1945592 and the Simons Foundation via the Simons Fellowship.}

\date{}

\begin{abstract}
We generalize the notion of an independent family on $\omega$ to $\mathcal{J}$-independent families on an arbitrary cardinal $\kappa$ and then show that under determinacy-like hypotheses that there can be no such families for reasonable $\mathcal{J}$. This is quite atypical as all other generalizations of the kinds of combinatorial families coming from the theory of cardinal characteristics follow a very different pattern, which will also be explained; in this sense independent families seem particularly special. 
The base theory assumed throughout is 
$$\textsf{ZF + }  \text{All infinite sets are Dedekind infinite}.$$
\end{abstract}

\maketitle

\setcounter{section}{-1}
\section{Introduction}

\setcounter{subsection}{-1}
\subsection{Some History}
When studying the continuum problem one defines several kinds of families of reals\footnote{A real is an infinite co-infinite subset of $\omega=\{0,1,2,\ldots\}$; $\R:=\{x\subseteq\omega:|x|=|\omega\bcl x|=\omega\}.$} which have either some kind of internal structure (almost disjoint families, independent families, and ultrafilter bases) and others which have some kind of maximality property (unbounded families, dominating families, splitting families, reaping families, ect.) one then associates to each kind of family a cardinal which is provably uncountable, though consistently strictly below the continuum. Here we will be focusing primarily on independent families and their closest relative, the reaping family. Fichtenholz and Kantorovich conceptualized the independent family in the early twentieth century \cite{FK}. Hausdorff then noticed shortly after that every cardinal $\kappa$ carries an independent family of cardinality $2^\kappa$ \cite{Haus}. A maximally independent family within a boolean algebra is essentially a maximal set of free generators of that algebra, however connections to set-theoretic topology have also been made.

Many of these connections are between independent families generally and various topological constructs, a connection of \textit{maximal} independent families is the notion of resolvability. A topological space is $\kappa$-resolvable if it can be partitioned into $\kappa$-many dense sets. Between \cite{app} and \cite{Hu} it is shown that if $\kappa$ is a regular cardinal and that there are no maximally independent families on $\kappa$ then for any $\gamma\ge\kappa$ any dense subspace of $\square_\kappa\ 2^\gamma$ of cardinality $\le\kappa$ is both $2$-resolvable and $\omega$-resolvable\footnote{The space $\square_\kappa 2^\gamma$ is the topology on $2^\gamma$ generated by basic open sets of the form $U_f = \{g\in 2^\gamma:g\supseteq f\}$ where $f\in \bigcup\{2^A:A\in[\gamma]^{<\kappa}\}$.}. Further, in \cite{app2} the authors use a particular maximally independent family to construct a dense subspace of $2^{2^\kappa}$ for which every open subset of this dense subspace has a dense subset of order $\kappa$, though no subset of this dense subspace is $2$-resolvable.

\subsection{Incorporating Determinacy}
For those families with internal structure it is unclear that there would be any maximal such families under the axiom of determinacy, in \cite{zach} Norwood and Neeman showed that there are not maximal almost disjoint on families on $\omega$ under Woodin's \thry{AD}$^+$ and in \cite{snt} Schrittesser and T\"ornquist sharpened this result showing that if $\thry{DC}_\R$ and \thry{Ramsey Uniformization} hold and $\omega\to(\omega)_2^\omega$ then there are no maximally almost disjoint families on $\omega$. 

One can then generalize the notions of all of these families to arbitrary cardinals. Chan, Jackson, and Trang generalized the notion of a maximal almost disjoint family to an arbitrary cardinal \cite{detad}.  In \cite{detad} the following results were obtained; note that all hypotheses,  except $V=L(\R)$, hold under Woodin's \thry{AD}$^+$.

\begin{thm*}[Chan, Jackson, Trang \cite{detad}] Assume that $\kappa$ is a cardinal. 
\begin{enumerate}
\item[\textup{(1)}]{If  $\omega<
cf(\kappa)<\Theta$ and the bold-face \thry{GCH} holds at $\cf(\kappa)$ then there are no maximal $\sB(\kappa)$-almost disjoint families and no maximal $[\kappa]^{<\kappa}$-almost disjoint families.}
\item[\textup{(2)}]{If $V=L(\R)+\thry{ AD}$ and $\cf(\kappa)\ge\Theta$ then there are well-orderable maximal $\sB(\kappa)$-almost disjoint families and well-orderable maximal $[\kappa]^{<\kappa}$-almost disjoint families. }
\item[\textup{(3)}]{ Assume $\thry{AD + DC}_\R$, $\omega\to(\omega)^\omega_2$, $\cf(\kappa)=\omega$, \thry{Ramsey Uniformization} holds for $\kappa$ and $\omega$ injects into every infinite subset of $\wp(\kappa)$, then there are no maximal $\sB(\kappa)$-almost disjoint families. }
\end{enumerate}
\end{thm*}

In this paper we generalize the notion of an independent family as follows:
\begin{defn} Let $\sJ$ be an ideal on an ordinal $\kappa$. We say $\sF\subseteq\sJ^+$ is $\sJ$-independent if for all finite disjoint $X,Y\subseteq\sF$
$$\thry{b}(X,Y):=\bigcap X\bcl\bigcup Y\in\sJ^+.$$ $\sF$ is a maximal $\sJ$-independent family if there is no $\sG\supsetneq \sF$ which is $\sJ$-independent.
\end{defn}
So that the classic independent family is a $\thry{fin}$-independent family on $\omega$. We obtain the following results.

\begin{cor*} \textbf{\ref{cort}.}
Assume $\thry{AD + DC}_\R$ and $\kappa<\Theta$ is a cardinal. Suppose $\sJ$ is an ideal on $\kappa$ such that either $\sJ$ is $\omega_1$-complete or $\cf(\kappa)=\omega$ and $\sB(\kappa)\subseteq\sJ\subseteq[\kappa]^{<\kappa}$.
There are no infinite maximal $\sJ$--independent families on $\kappa$.\\
\end{cor*}

\begin{cor*} \textbf{\ref{corm}.}
Assume $J$ is a set of ordinals and $V=L(\R,J)+\thry{AD}+\thry{DC}_\R$ and $\kappa$ is a cardinal. Suppose $\sJ$ is an ideal on $\kappa$ such that either $\sJ$ is $\omega_1$-complete or $\cf(\kappa)=\omega$ and $\sB(\kappa)\subseteq\sJ\subseteq[\kappa]^{<\kappa}$.
There are no infinite maximal $\sJ$--independent families on $\kappa$.
\end{cor*}

Then one can use these corollaries to establish that the most general $V=L(\wp(\R))\thry{ + AD}^+$ models do not have maximally independent families. 

\begin{thm*} \textbf{\ref{lpr}.} 
Suppose $V=L(\wp(\R))\thry{ + AD}^+$ and $\sJ$ is an ideal on a cardinal $\kappa$ such that either $\sJ$ is $\omega_1$-complete or $\cf(\kappa)=\omega$ and $\sB(\kappa)\subseteq\sJ\subseteq[\kappa]^{<\kappa}$.
If $\sF$ is an infinite $\sJ$--independent family on $\kappa$, then $\sF$ is not maximal. 
\end{thm*}

\section{Determinacy and $L(\wp(\R))$}

In this section we merely recount the necessary facts about determinacy to establish our results.

\begin{fact}
Assume \thry{AD} then there are no injections $\omega_1\to\R$, for every $A\subseteq\R$ there is an open set $U$ such that $A\triangle U$ is meagre, and the class of meagre sets is closed under well-ordered unions.
\end{fact}

\begin{defn}
$\Theta$ is the least ordinal onto which one cannot surject the reals. 
\end{defn}

The following lemma, due to Moschovakis, will be used to conclude some of the main corollaries. Theorem \ref{main} has stronger hypotheses than what follows from \thry{AD} alone, but using this lemma allows the theorem to apply to $\kappa<\Theta$. 

\begin{lem}[Moschavakis \cite{mos}]
Assume \thry{AD} and $\kappa<\Theta$. $\R$ surjects onto $\wp(\kappa)$.
\end{lem}

Assuming $\thry{AD}^+$ we understand well-ordered subsets of $\wp(\kappa)$ for $\kappa<\Theta$.

\begin{defn}
Suppose $\kappa$ is a cardinal. The bold-face \thry{GCH} at $\kappa$ holds if there is no injection $\kappa^+\to\wp(\kappa)$. 
\end{defn}

\begin{fact}[Steel, Woodin \cite{steel}]
Suppose \thry{AD}$^+$ and $\kappa<\Theta$ is a cardinal. The bold-face \thry{GCH} holds at $\kappa.$
\end{fact}

We will analyze models of the form $V=L(\wp(\R))+\thry{AD}^+$. The following theorem, due to Woodin, shows that such models come in two forms.

\begin{thm}[Woodin]
Suppose $V=L(\wp(\R))+\thry{AD}^+$ then either $\thry{AD}_\R$ holds or there is a set of ordinals $J$ such that $V=L(\R,J)+\thry{AD}^+.$
\end{thm}

Considering the above and after analyzing models of the form $L(\R)$, to analyze the $\thry{AD}^+$ models of the form $L(\wp(\R))$ we will only need to consider the case that $\thry{AD}_\R \thry{ + }\neg\thry{DC}$ holds. In such a case it is helpful to know the cofinality of $\Theta$.

\begin{lem}[Solovay \cite{solovay}]
If $V=L(\wp(\R))\thry{ + AD}^+\thry{ + AD}_\R\thry{ + }\neg\thry{DC}$, then $\cf(\Theta)=\omega$.
\end{lem}

A hypothesis of theorem \ref{main} is $\thry{DC}$, which does not follow generally from $\thry{AD}^+$. To work around this we apply \ref{main} within inner models of $\thry{AD}^+$ where $\thry{DC}_\R$ will suffice. 

\begin{prop}\label{HODDC}
Let $J$ be a set of ordinals and $C\subseteq\R$, then  $L(\R,J)\model\thry{DC}_\R\to\thry{DC}$ and $\thry{HOD}_{\R\cup\{C,J\}}\model\thry{DC}_\R\to\thry{DC}$.  
\end{prop}
\begin{proof}
The argument is essentially the same for either $L(\R,J)$ or $\thry{HOD}_{\R\cup\{C,J\}}$. Let $\sM$ be either $L(\R,J)$ or $\thry{HOD}_{\R\cup\{C,J\}}$.

Let $X\neq 0$ be a nonempty set in $\sM$. There is a surjection $\pi:\R\times\thry{On}\twoheadrightarrow X$ in $\sM$. Let $R$ be a relation on $X$ with $\dom R=X$. Define a $\sigma:\R\times\thry{On}\to\thry{On}$ by 
$$\sigma(a,\alpha) = \min\{\beta:(\exists b\in \R)( \langle \pi(a,\alpha),\pi(b,\beta)\rangle\in R)\}$$ and a tree $T\subseteq\R^{<\omega}$ with a function $\tau:T\to \thry{On}$ recursively. Let $x\in X$ be given and $a_0\in\R$ and $\alpha_0\in \thry{On}$ such that $\pi(a_0,\alpha_0)=x$. Then $T_1=\{\langle a_0\rangle\}$ and $\tau(\langle a_0\rangle)=\alpha_0$.

Once $T_n\subseteq \R^n$ has been defined say 
$\langle a_0,\ldots,a_{n-1},a_n \rangle \in T_{n+1}$ if and only if both $\langle a_0,\ldots,a_{n-1}\rangle\in T_{n}$ and, letting $\beta = \tau(\langle  a_0,\ldots,a_{n-1}\rangle)$, $$\langle \pi(a_{n-1},\beta), \pi(a_n, \sigma(a_{n-1},\beta))\rangle\in R$$ and let 
$\tau(\langle a_0,\ldots,a_n\rangle) = \sigma(a_{n},\sigma(a_{n-1},\beta)).$ Let $T=\bigcup\{T_n:n<\omega\}$. Notice that $T$ has no terminal nodes since for any $\langle a_0,\ldots,a_{n}\rangle\in T$ there is some $b\in\R$ such that $$\langle \pi(a_n,\sigma(a_{n-1},\tau(\langle a_0,\ldots,a_{n-1}\rangle )), \ \pi(b, \sigma(a_n,\sigma(a_{n-1},\tau(\langle a_0,\ldots, a_{n-1}\rangle))))\rangle\in R$$ and any such $b$ has $\langle a_0,\ldots,a_n,b\rangle\in T_{n+1}.$ So by \thry{DC}$_\R$ there is a $\langle a_n:n<\omega\rangle\in\R^\omega$ with $\langle a_n:n<\omega\rangle\in[T]$. Let $\alpha_n = \tau(\langle a_0,\ldots,a_n\rangle)$ for all $n<\omega$ and let $\langle x_n:n<\omega\rangle$ be defined by 
$x_n = \pi(a_n, \alpha_n).$ We show that $\langle x_n,x_{n+1}\rangle\in R$ for all $n$ and $x_0=x$. Well, $\tau(\langle a_0\rangle)=\alpha_0$ and so $x_0 =  \pi(a_0,\alpha_0)=x.$ Now, notice that 
$\alpha_{n} =\sigma(a_{n-1},\alpha_{n-1})$
and since $\langle a_0,\ldots,a_{n}\rangle\in T$, 
$$\langle x_{n-1},x_n\rangle =  \langle \pi(a_{n-1},\alpha_{n-1}),\ \pi(a_n, \sigma(a_{n-1},\alpha_{n-1}))\rangle\in R.$$
\end{proof}

Then to show that there are well-orderable reaping families on suitable $\kappa$ in $L(\R,J),$ we need to know how $L(\R,J)$ interacts with its \thry{HOD}. 

\begin{thm}[Woodin] \label{woodin}
Suppose $V=L(\R,J)+\thry{AD}$. There is a forcing $\A\in\thry{HOD}_J$ such that $\thry{HOD}_J\model |\A|=\Theta$ and such that for all $S\subseteq\thry{On}$ there is a $g\subseteq\A$ which is generic over $\thry{HOD}_J$ and $S\in \thry{HOD}_J[g]$. 
\end{thm}

In $L(\R,J)$ every set is definable from a real and $J$. $\A$ is essentially the Vop$\check{\text{e}}$nka algebra for adding a real, the forcing consisting \thry{OD}$_J$ $\infty$-Borel codes for sets of reals. Chan in \cite{detcomb} has detailed accounts of this theorem. For us, the only facts about $\A$ we will use is that each $p\in\A$ is a code for a set of reals so that in $L(\R,J)$, $p$ is essentially a set of reals which we will not distinguish from its code in $\thry{HOD}_J$.

\section{Reaping Families}

In this section we focus on reaping families due to their connections with independent families, though the analysis is the same with the same tools for any of: reaping families, dominating families, unbounded families, splitting families, almost disjoint families, towers, families with the strong finite intersection property but no psuedo-intersection, and shattering families. 

\begin{defn}
Suppose $\sJ$ is an ideal on an ordinal $\kappa$. Then $y\subseteq\kappa$ $\sJ$-splits $x\in\sJ^+$ if both $x\cap y\in\sJ^+$ and $x\bcl y\in\sJ^+$. Then $y\subseteq\kappa$ $\sJ$-splits $\sR\subseteq\sJ^+$ if $y$ $\sJ$-splits every $x\in \sR$. We say $\sR$ is $\sJ$-reaping if there is no single $y\subseteq\kappa$ which $\sJ$-splits $\sR$. 
\end{defn}

It is a classic theorem that for $\kappa=\omega$, $\sJ=\thry{fin}$ that every countable $\sR\subseteq[\omega]^\omega$ is split, and thus all reaping families on $\omega$ are uncountable. The argument would enumerate any countable subset of $[\omega]^\omega$ in order-type $\omega$, then go through this enumeration to construct a set which splits the entire sequence, picking integers to put into the splitting set and others to keep out. This argument will be adapted to show that there are no well-ordered reaping families on certain $\kappa$. 

\begin{defn}
For a cardinal $\kappa$, we define two ideals which coincide should $\kappa$ be regular.
$\sB(\kappa)=\bigcup\{\wp(\eta):\eta<\kappa\}$ and 
$[\kappa]^{<\kappa}=\{x\subseteq\kappa:|x|<\kappa\}.$
\end{defn}

\begin{lem}
Suppose $\kappa$ is a regular cardinal and the bold-face \thry{GCH} holds at $\kappa$ then there are no well-ordered $\sB(\kappa)$-reaping families on $\kappa$. 
\end{lem}
\begin{proof}
Suppose $\sR\subseteq\sB(\kappa)^+$ is well-ordered, we split it. Since $\sR$ is well-orderable there is an enumeration of $\sR$ in order-type $\kappa$, $\langle x_\xi:\xi<\kappa\rangle$. Define sequences $\langle A_\xi:\xi<\kappa\rangle$ and $\langle B_\xi:\xi<\kappa\rangle$ of bounded subsets of $\kappa$ such that for all $\zeta<\xi<\kappa$
\begin{equation}
\sup A_\zeta <\min B_\zeta<\sup B_\zeta <\min A_\xi\tag{$\star$}
\end{equation}
and $A_\xi\cap x_\zeta\neq 0\neq B_\xi\cap x_\zeta.$ Once $\langle A_\xi:\xi<\eta\rangle$ and $\langle B_\xi:\xi<\eta\rangle$ have been defined obeying $(\star)$. Let $\vartheta = \sup\{\sup(B_\xi):\xi<\eta\}+1,$ since $\eta<\kappa$ and $\kappa$ is regular $\vartheta<\kappa$. Define
$A_\eta = \{\min(x_\xi\bcl\vartheta):\xi<\eta\}$ again, since $\kappa$ is regular, 
$\delta:=\sup A_\eta + 1$ has $\delta<\kappa$ and let 
$B_\eta = \{\min(x_\xi\bcl\delta):\xi<\eta\}.$

Then $y=\bigcup \{A_\xi:\xi<\kappa\}$ splits $\sR$. To see this, for any $\eta<\kappa$ and any $\xi>\eta$, $A_\xi\cap x_\eta\neq 0$ and in fact the map $\xi\mapsto\min(x_\eta\cap A_\xi)$ is an injection and so $\kappa\bcl\eta$ injects into $y\cap x_\eta$, the same can be said of $\bigcup\{B_\xi:\xi<\kappa\}$ and $$\bigcup\{A_\xi:\xi<\kappa\}\cap\bigcup\{B_\xi:\xi<\kappa\}=0.$$ 
\end{proof}

\begin{thm}\label{noreaping}
Suppose $\kappa$ is a cardinal and the bold-face \thry{GCH} holds at $\cf(\kappa)$. There are no well-orderable $\sB(\kappa)$-reaping families on $\kappa$ and if in addition $\cf(\kappa)>\omega$ then there are no well-orderable $[\kappa]^{<\kappa}$-reaping families on $\kappa$. 
\end{thm}
\begin{proof}
First we assume $\cf(\kappa)>\omega$ and show that there are no well-orderable $[\kappa]^{<\kappa}$-reaping families on $\kappa$. Suppose $\sR\subseteq[\kappa]^\kappa$ was well-orderable, we construct a splitting set. Fix an increasing continuous cofinal $\rho:\cf(\kappa)\to\kappa$. For $x\in\sR$ let 
$$\pi(x)=\{\xi\in \cf(\kappa) : \cf(\rho(\xi))=\omega\ \land (\forall \eta<\xi)( \ |x\cap\rho(\xi)\bcl \rho(\eta)|=|\rho(\xi)|)\}.$$ It is perhaps less clear that for all $x\in \sR$ that $\pi(x)\in\sB(\cf(\kappa))^+$. Let $x\in\sR$ and $\eta_0<\cf(\kappa)$ be given. We construct an increasing sequence. Once $\eta_i$ is defined, since $|x|=\kappa$, there is some least $\eta_{i+1}>\eta_i$ such that $|x\cap \rho(\eta_{i+1})\bcl\rho(\eta_i)|\ge |\rho(\eta_i)|$. Let $\eta=\sup\{\eta_i:i<\omega\}$, then $\cf(\eta)=\omega$ and so $\rho(\eta) = \sup\{\rho(\eta_i):i<\omega\}$ is also of countable cofinality as $\rho$ is continuous. Further for any $i<j<\omega$, $|x\cap \rho(\eta)\bcl\rho(\eta_i)|\ge|\rho(\eta_j)|$ and so for all $i<\omega$, $|x\cap \rho(\eta)\bcl\rho(\eta_i)|=|\rho(\eta)|.$ So, we have found a member of $\pi(x)$ larger than any given $\eta_0$. 

Since $\sR$ is well-orderable, so is $\pi"\sR$ and so we enumerate $\pi"\sR$ in order-type $\cf(\kappa)$ as $\langle x_\eta:\eta<\cf(\kappa)\rangle$ and as in the proof of the previous we get sequences of bounded subsets of $\cf(\kappa)$, $\langle A_\eta:\eta<\cf(\kappa)\rangle$ and $\langle B_\eta:\eta<\cf(\kappa)\rangle$ such that for all $\zeta<\eta<\cf(\kappa)$ we have $x_\zeta\cap A_\eta\neq 0\neq x_\zeta\cap B_\eta$ and 
$\sup A_\zeta+1<\min B_{\zeta}<\sup B_{\zeta}+1<\min A_\eta.$

For each $\xi<\cf(\kappa)$ let $$I_\xi=[\rho(\min A_\xi), \rho(\sup A_\xi+1)]\text{ and }J_\xi=[\rho(\min B_\xi),\rho(\sup B_{\xi}+1)].$$
Notice this, if $x\in \sR$ and $\pi(x)=x_\eta$ and $\eta<\xi$ then there is some $\alpha\in A_\xi\cap \pi(x)$ so that $|x\cap\rho(\alpha)\bcl \rho(\min A_\xi)|=|\rho(\alpha)|.$ In particular, whenever $x\in\sR$ and $\pi(x)=x_\eta$, for all $\xi>\eta$, $|x\cap I_\xi|\ge |\rho(\min A_\xi)|$ and since $\langle \min(A_\xi):\xi<\cf(\kappa)\rangle$ runs cofinally through $\cf(\kappa)$, for any $x\in \sR$, 
$|x\cap\bigcup\{I_\xi:\xi<\cf(\kappa)\}|=\kappa.$ The same can be said of the $J_\xi$'s, so letting $y=\bigcup\{I_\xi:\xi<\cf(\kappa)\}$ we have that $y$ $[\kappa]^{<\kappa}$-splits $\sR$ as $\{J_\xi:\xi<\cf(\kappa)\}\cap\bigcup\{I_\xi:\xi<\cf(\kappa)\}=0$.  

Showing that there are no well-orderable $\sB(\kappa)$-reaping families on $\kappa$ is similar. Let $\pi(x)=\{\xi\in\cf(\kappa):x\cap[\rho(\xi),\rho(\xi+1))\neq 0\}$ instead and if $\hat{y}\subseteq\cf(\kappa)$ splits $\pi"\sR$ then $y:=\bigcup\{[\rho(\eta),\rho(\eta+1)):\eta\in\hat{y}\}$ splits $\sR$. 
\end{proof}

\begin{cor}
Assume $\thry{AD}^+$ and $\kappa$ is a cardinal with $\cf(\kappa)<\Theta$. There are no well-orderable $\sB(\kappa)$-reaping families on $\kappa$. If $\cf(\kappa)>\omega$ then there are no well-orderable $[\kappa]^{<\kappa}$-reaping families on $\kappa$.
\end{cor}

\begin{qstin}
Can there be well-orderable $[\kappa]^{<\kappa}$-reaping families on a cardinal $\kappa$ of countable cofinality under $\thry{AD}^+$?
\end{qstin}

Similar arguments apply to unbounded, dominating, and splitting families, among other kinds of generalizations of families coming from the theory of cardinal characteristics. Bold-face \thry{GCH} at $\cf(\kappa)$ is sufficient to prove that there are no well-orderable such families on $\kappa$ by pushing down along a cofinal sequence, using constructions at the level of $\cf(\kappa),$ and then push back up along a cofinal sequence. The constructions here at the level of $\cf(\kappa)$ are similar to the classic constructions showing that there are no countable such families on $\omega$. 

In the \thry{AD} context there is very little that can be determined about the structure of $\wp(\kappa)$ if $\cf(\kappa)\ge\Theta$. So, to analyze the structure of $\wp(\kappa)$ for such cardinals model-dependent analysis should be done and the most classic models of \thry{AD} are models of the form $L(\R,J)$ for some set of ordinals $J$. Here we see that there will be well-orderable reaping, dominating, unbounded, and splitting families. The core tool for this is Woodin's theorem \ref{woodin} that such models are symmetric collapse extensions of their \thry{HOD}. See Chan's \cite{detcomb} for a detailed exposition. In \cite{detad} Chan, Jackson, and Trang use these facts to establish that if $V=L(\R)+\thry{AD}$ and $\cf(\kappa)\ge \Theta$ then there are well-orderable maximal $\sB(\kappa)$-almost disjoint families and well-orderable $[\kappa]^{<\kappa}$-almost disjoint families, both well-orderable in order-type $\kappa^+$. Their argument is that in \thry{HOD}$^{L(\R)}$ there are such families which are maximal, they are then well-ordered in \thry{HOD}$^{L(\R)}$ in order-type $(\kappa^+)^{\thry{HOD}}=\kappa^+$. Then one notices that the maximality of such families are indestructible by forcings of cardinality $<\cf(\kappa)$ and so remain maximal in $V=L(\R)$. The same argument can be applied to any other of the families, except independent families. Unbounded, dominating, and splitting families also enjoy small forcing indestructibility. While reaping families do not in general enjoy small-forcing indestructibility, in any model of \thry{ZF} and any $\kappa$ in this model there is a reaping $\sR\subseteq\wp(\kappa)$ which cannot be split by small forcing. 

\begin{prop}\label{yesreaping}
Suppose that $V=L(\R,J)+\thry{AD}^+$ and $\kappa$ is a cardinal such that $\cf(\kappa)\ge\Theta$ and $\sJ=\sB(\kappa)$ or $\sJ=[\kappa]^{<\kappa}$. There is a well-orderable $\sJ$-reaping family on $\kappa$ of cardinality $2^{\kappa,\thry{HOD}_J}.$
\end{prop}
\begin{proof}
Notice that $\sJ^{\thry{HOD}_J}= \sJ\cap\thry{HOD}_J$. Let $$\sR=\{a\subseteq\kappa:a\notin\sJ\land \kappa\bcl a\notin\sJ\}\cap\thry{HOD}_J\in\thry{HOD}_J$$ So, $\thry{HOD}_J\model|\sR|=2^\kappa.$
Then $\sR$ is well-ordered in order-type $(2^\kappa)^{\thry{HOD}_J}$ and remains unsplit in $V=L(\R,J).$ Let $\A$ denote the Vop$\check{\text{e}}$nka algebra for adding a real. To see that $\sR$ remains unsplit it suffices to show that whenever $\dot{x}$ is an $\A$-name for a subset of $\kappa$ and $p\force \dot{x}\in\dot{\sJ}^+$ then there is a $q\le p$ such that $\{\eta<\kappa:q\force \check{\eta}\in \dot{x}\}\in\sJ^+.$ For then any $x\in\sJ^+$ in $V$ has a subset $y\in\sR$. It suffices to show that $\A$ has the $\Theta$-c.c.

Working in $V$, each $p\in\A$ is a non-empty set of reals, see \cite{detcomb}. If $\langle p_\xi:\xi<\Theta\rangle$ in $\thry{HOD}_J$ is a maximal antichain in $\A$ then in $V$, $\bigcup\{p_\xi:\xi<\Theta\}=\R$ for if this union is nonempty then $\R\bcl\{p_\xi:\xi<\Theta\}$ is a condition in $\A$ and is incompatible with each $p_\xi$. It is also the case that since $\{p_\xi:\xi<\Theta\}$ forms an antichain, for $\xi<\eta<\Theta$,  $p_\xi\cap p_\eta=0$ in $V,$ when regarded as sets of reals. Then in $V$ we define $\phi:\R\to\Theta$ by $\phi(r)$ is the $\xi$ such that $r\in p_\xi$, then $\phi$ is a surjection of $\R$ onto $\Theta$, thus $\Theta<\Theta$.
\end{proof}

We will see in the next section that it is not the case that all reaping families are indestructible by small forcings, examples of reaping families which can be split by small forcing extensions will naturally arise from independent families.

\begin{qstin}
Assume $V=L(\thry{On}^\omega)$\thry{ + AD}$^+$ or some other Chang-like model. Can there be maximal almost disjoint families? Can there be well-orderable dominating, unbounded, splitting, and reaping families?
\end{qstin}

In the Chang-like models the algebra to add a subset of $\kappa$ to \thry{HOD} will generally not be of cardinality $<\cf(\kappa)$ so small forcing indestructibility does not suffice.  

\section{Independent Families}

\begin{defn} For a $\sJ$-independent family $\sF$ let $\sR(\sF)$ consist of the set of boolean combinations of members of $\sF$:
$$\sR(\sF):=\{\thry{b}(X,Y):X,Y\subseteq \sF\text{ finite disjoint}\}$$
\end{defn}

The connection to reaping families is the same motivation for the notation. 

\begin{prop}
Let $\sF$ be an $\sJ$-independent family. Then $\sF$ is maximally $\sJ$-independent if and only if $\sR(\sF)$ is $\sJ$-reaping
\end{prop}
\begin{proof}
For any $x\in\sJ^+$, $\sF\cup\{x\}$ is $\sJ$-independent if and only if $x$ $\sJ$-splits $\sR(\sF)$. 
\end{proof}

Since any well-orderable $\sJ$-independent family $\sF$ has that $\sR(\sF)$ is well-orderable, theorem \ref{noreaping} immediately gives the following. 

\begin{prop}
Assume the bold-face \thry{GCH} holds at $\cf(\kappa)$ then there are no maximal well-orderable $\sB(\kappa)$-independent families on $\kappa$ and if $\cf(\kappa)>\omega$ nor are there any maximal well-orderable $[\kappa]^{<\kappa}$-independent families on $\kappa$.
\end{prop}

One might object that the generalization here is the incorrect one. In the classic case, an independent family on $\omega$ is a $\thry{fin}$-independent family where boolean combinations are of length $<\omega$. One might then insist that when generalizing the notion to arbitrary cardinals $\kappa$ we should be allowed $<\lambda$ length boolean combinations for some $\omega<\lambda\le\kappa$.\footnote{In \cite{app} and \cite{Hu} the notion of independence used is one in which $\sF\subseteq\wp(\kappa)$ is said to be independent if for all $X,Y\subseteq\sF$ disjoint with $|X\cup Y|<\kappa$, $\bigcap X\bcl\bigcup Y$ is infinite.} Under \thry{AD} the existence of such families is immediately refuted.

\begin{prop} 
Assume $\sF\subseteq\wp(\kappa)$ is an infinite family of subsets of $\kappa\in\thry{On}$ such that whenever $X,Y\in[\sF]^{\le\omega}$ are disjoint $\thry{b}(X,Y)\neq 0$. There is a well-ordering of the reals in order-type $\le\kappa$. \footnote{We assume that all infinite sets are Dedekind infinite only to avoid writing ``Assume there is an injection $\omega\to \sF$'', which we would write every time otherwise.}
\end{prop}
\begin{proof}
Let $\langle x_n:n<\omega\rangle$ be an injection $\omega\to\sF$. For each $f:\omega\to 2$ let 
$$A^f = \bigcap \{x_n:f(n)=0\}\bcl\bigcup\{x_n:f(n)=1\},$$ the hypothesis that $\omega$-length boolean combinations are non-empty gives that for all $f\in 2^\omega$ that $A^f\neq 0$. Moreover, if $f,g\in 2^\omega$ are distinct, then $A^f\cap A^g=0$ and so $\Phi:2^\omega\to\kappa$ by $\Phi(f)=\min(A^f)$ is an injection. 
\end{proof}

Many of our results showing that there are no maximally independent families in determinacy contexts really show that there are no infinite maximally independent families. Let us briefly address finite independent families. 

\begin{defn}
Say an ideal $\sJ$ on a set $X$ is non-atomic if whenever $b\in\sJ^+$ then there is an $a\subsetneq b$ with both $a\in\sJ^+$ and $b\bcl a\in \sJ^+$.
\end{defn}

\begin{prop}
Let $\sJ$ be a non-atomic ideal on a set $X$ and $\sF$ a finite $\sJ$-independent family on $X$. $\sF$ is not $\sJ$-maximal. 
\end{prop}
\begin{proof}
Say $\langle x_i:i<n\rangle$ is an injective enumeration of $\sF$. For each $s\in 2^n$ let 
$$B_s = \bigcap\{x_i:s(i)=0\}\bcl\bigcup\{x_i:s(i)=1\}\in\sJ^+$$ Notice that if $s\neq t$ then $B_s\cap B_t=0$. For each $s\in 2^n$ let $a_s\subsetneq B_s$ be such that both $a_s\in\sJ^+$ and $B_s\bcl a_s\in\sJ^+$. Let $y=\bigcup\{a_s:s\in 2^n\}$. Then $\sF\cup\{y\}$ is $\sJ$-independent. Let $B\in\sR(\sF)$ be given and $s\in 2^n$ be such that $B_s\subseteq B$, then $B\cap y\supseteq B_s\cap y = a_s\in \sJ^+$ and 
$B\bcl y\supseteq B_s\bcl y = B_s\bcl a_s\in\sJ^+.$
\end{proof}

A key tool in showing that under determinacy-like hypotheses that there are no maximally independent families on an ordinal is the density ideal associated to a family. The density ideal goes back to Shelah in his paper showing the consistency of $\thry{Con}(\mathfrak{i}<\mathfrak{u})$ \cite{shelah} but in \cite{higher} Fischer and Montoya generalize the notion to an $\sB(\kappa)$-independent family on a regular cardinal $\kappa$ as considered here. 

\begin{defn}
If $\sF\subseteq\wp(\kappa)$ is $\sJ$-independent then the $\sJ$-density ideal of $\sF$ is 
$$\id_\sJ(\sF) = \{x\subseteq\kappa : (\forall B\in\sR(\sF))(\exists C\in \sR(\sF))(C\subseteq B\land C\cap x\in\sJ)\}.$$
\end{defn}

\begin{fact}
If $\sF$ is a $\sJ$-independent family then $\id_\sJ(\sF)\supseteq\sJ$ and $\sR(\sF)\cap\id_\sJ(\sF)=0$. 
\end{fact}

\begin{prop}[Fischer, Montoya \cite{higher}] If $\sF\subseteq\sG$ are $\sJ$-independent families then $\id_{\sJ}(\sF)\subseteq\id_{\sJ}(\sG).$
\end{prop}
\begin{proof}
Suppose $D\in\id_{\sJ}(\sF)$. Let $B\in \sR(\sG)$ be given and $\hat{X},\hat{Y}\subseteq\sG$ be finite disjoint such that $\thry{b}(\hat{X},\hat{Y})=B$. Let $X=\hat{X}\cap \sF$ and $Y=\hat{Y}\cap\sF$. Then since $D\in\id_{\sJ}(\sF)$ there are $X',Y'\subseteq\sF$ finite such that $(X\cup X')\cap (Y\cup Y')=0$ and $\thry{b}(X\cup X', Y\cup Y')\cap D\in \sJ.$ Then $C:=\thry{b}(\hat{X}\cup X',\hat{Y}\cup Y')\in\sR(\sG)$ has $C\subseteq B$ and 
$$C\cap D\subseteq \thry{b}(X\cup X', Y\cup Y')\cap D\in \sJ.$$
\end{proof}

\begin{lem}
Suppose there are no injections $\omega_1\to\R$. Let $\sJ$ be any ideal on an ordinal $\kappa$ and $\sF$ be an infinite $\sJ$-independent family on $\kappa$. There is a partition $\langle A_n:n\in\omega\rangle$ of $\kappa$ with $A_n\in\id_\sJ(\sF)$ for all $n\in\omega.$ 
\end{lem}
\begin{proof}
Fix an injection $\langle x_n:n\in\omega\rangle$ of $\omega\to\sF$. For each $f\in 2^\omega$ let $$A^f = \bigcap \{x_n:f(n)=0\}\bcl\bigcup\{x_n:f(n)=1\}.$$ We show that for all $f\in 2^\omega$ that $A^f\in\id_\sJ(\sF)$ and to show this, it suffices to show that $A^f\in\id_{\sJ}(\{x_n:n\in\omega\})$. For any given $B\in\sR(\{x_n:n<\omega\})$ and $f\in 2^\omega$, there are unique finite disjoint $X,Y\subseteq\{x_n:n\in\omega\}$ such that $B=\thry{b}(X,Y)$. Let $k$ be sufficiently large such that $X\cup Y\subseteq\{x_n:n<k\}.$
If $f(k)=0$ let $C=\thry{b}(X,Y\cup\{x_k\})$ and if $f(k)=1$ let $C=\thry{b}(X\cup\{x_k\},Y)$, in either case we have that $C\subseteq B$ and $A^f\cap C=0$. It is not hard to see that $\langle A^f:f\in 2^\omega\rangle$ partitions $\kappa$.

To finish the proof we show that $E=\{f\in 2^\omega: A^f\neq 0\}$ is countable. Inject $E$ into $\kappa$ by $f\mapsto\min(A^f)$, then by the hypothesis that there is no injection $\omega_1\to\R$ we must have that the order-type of the range of $f\mapsto\min(A^f)$ is countable. 
\end{proof}

\begin{remark} 
Notice that the above really shows that, in \thry{ZF}, if $\sF$ is an independent family into which $\omega$ injects, then there is a partition $\langle A^x:x\in\R\rangle$ of $\kappa$ with each $A^x\in\id_\sJ(\sF)$. The hypothesis that $\omega_1\not\le|\R|$ gives that this partition is countable.
\end{remark}

We can already see that independent families stand out. Assuming $V=L(\R)$, all the other kinds of the combinatorial families coming from the theory of cardinal characteristics exist on cardinals $\kappa$ with $\cf(\kappa)\ge\Theta$, we even have well-orderable such families. However, the following shows that this argument that applies to most other kinds of these families, small forcing indestructibility together with theorem \ref{woodin}, has no chance of being adapted to independent families. 

\begin{thm}\label{nowo}Assume $\omega_1\not\le|\R|$. 
Suppose $\sJ$ is an $\omega_1$-complete ideal on $\kappa\in\thry{On}$ and $\sF$ is an infinite well-orderable $\sJ$-independent family. $\sF$ is not $\sJ$-maximal. 
\end{thm}
\begin{proof}
We split $\sR(\sF)$. Begin with a partition $\langle A_n:n\in\omega\rangle$ of $\kappa$ with $A_n\in\id_{\sJ}(\sF)$ for all $n\in\omega$. Define 
$\pi:\sR(\sF)\to [\omega]^\omega$ by 
$\pi(B) = \{n\in\omega  : A_n\cap B\in\sJ^+\}.$ See that $\pi(B)$ is infinite for all $B\in\sR(\sF).$ For if $B\in\sR(\sF)$ is otherwise then there is $N\in\omega$ with $A_n\cap B\in\sJ$ for all $n\ge N$ and so by the $\omega_1$-completeness of $\sJ$ $$B\cap\bigcup\{A_n:n\ge N\}=\bigcup\{B\cap A_n : n\ge N\}\in\sJ\subseteq\id_\sJ(\sF).$$
but each $A_n\in\id_\sJ(\sF)$ and hence 
$B\cap\bigcup\{A_n:n<N\}\in\id_\sJ(\sF)$ and so 
$$B = \left(B\cap\bigcup\{A_n:n<N\} \right) \cup\left(B\cap\bigcup\{A_n:n\ge N\}\right)\in \id_\sJ(\sF)$$
which is absurd in light of the fact that $\sR(\sF)\cap\id_\sJ(\sF)=0$. 

Then since $\sF$ is well-orderable $\sR(\sF)$ is also well-orderable. Then $\pi$ is a surjection onto its range and so $\pi"\sR(\sF)$ is well-orderable. Since there are no injections $\omega_1\to\R$, it must be that $\pi"\sR(\sF)$ is countable. Let $\langle x_n:n<\omega\rangle$ enumerate $\pi"\sR(\sF)$. Then define recursively a strictly increasing $y:\omega\to\omega$ with $y(2n), y(2n+1)\in x_n$ for all $n\in\omega$. Let $z=\bigcup\{A_{y(2n)}:n\in\omega\}.$ The claim is that $z$ splits $\sR(\sF)$ so that $\sF\cup\{z\}$ is $\sJ$-independent. To see this, let $B\in\sR(\sF)$, say $\pi(B)=x_n$, then 
$$B\cap z\supseteq B\cap A_{y(2n)}\in\sJ^+\text{ and } B\bcl z\supseteq B\cap A_{y(2n+1)}\in\sJ^+.$$
\end{proof}

\begin{cor}[\thry{ZFC}] Let $\sF$ be an infinite $\sJ$-independent family on $\kappa$ for an ideal $\sJ$ which is $\cc^+$-complete. Suppose $\dot{\Pn}$ is a $\thry{Col}(\omega,\R)$-name for a notion of forcing and $\thry{Col}(\omega,\R)\force\dot{\Pn}\text{ adds a splitting real}.$ Then 
$\thry{Col}(\omega,\R)*\dot\Pn\force \check{\sF}\text{ is not }\dot{\sJ}\text{-maximal}.$ More generally, let $\Q$ be a notion of forcing which adds a splitting set of reals in that there is a $\Q$-name $\dot{E}$ such that for all infinite $X\in\wp^V(\R)$, 
$\Q\force \check{X}\cap\dot{E}\neq 0\neq \check{X}\bcl \dot{E},$ then 
$\Q\force\check{\sF}\text{ is not }\dot{\sJ}\text{-maximal}$
\end{cor}
\begin{proof}
In $V$ there is a partition $\langle A^x:x\in\R^V\rangle$ of $\kappa$ with each $A^x\in \id_\sJ^V(\sF)$. Define $\pi:\sR(\sF)\to\wp(\R)$ by $\pi(B) = \{x:A^x\cap B\in\sJ^+\}$. The completeness of $\sJ$ implies that $\pi(B)$ is infinite for all $B\in\sR(\sF)$. Let $g\subseteq\Q$ be generic and $E\subseteq\R$ in $V[g]$ be such that for all infinite $X\in\wp^V(\R),$ $X\cap E\neq0\neq  X\bcl E$. Let $z=\bigcup\{A^x : x\in E\}.$ Then $\sF\cup\{z\}$ is independent: given $B\in\sR(\sF)$ fix $x_0\in \pi(B)\cap E$ and $x_1\in \pi(B)\bcl E$, then 
$$B\cap z\supseteq B\cap A^{x_0}\in \sJ^+\text{ and } B\bcl z\supseteq B\cap A^{x_1}\in\sJ^+.$$
\end{proof}

\begin{remark} 
This corollary delivers on the previous promise to give examples of reaping families which are split by small forcings. If $\kappa$ is a cardinal with $\cf(\kappa)>\cc$ then for any maximal $\sB(\kappa)$-independent family $\sF$, $\sR(\sF)$ is split in any generic extension by $\thry{Col}(\omega,2^\cc)$, $\thry{Col}(\omega,\R)*\thry{Col}(\omega,\R),$ $\thry{Col}(\omega,\R)*\thry{Add}(\omega,1),$ or $\thry{Add}(\R,1)$. 
\end{remark}

So, assuming \thry{AD}$^+$ the only hope for a reasonable independent family is a non-well-orderable family. We work toward establishing that there can be no such thing on a cardinal $\kappa<\Theta$ and that in the canonical models of \thry{AD} there are none at all. 

\begin{lem}[$\Delta$-System Lemma] \label{delta}
Suppose \thry{ZF + AC}$_\omega(\R)$. Let $B\subseteq\R\times V$ be a relation such that for all $t\in\R$, $$B_t:=\{x:\langle t,x\rangle\in B\}$$ is finite. There is a $r\subseteq\ran B$ such that for all $n\in\omega\bcl 2,$ there is a $D\subseteq \R$ with $|D|>n$ and such that $\{B_t:t\in D\}$ forms a $\Delta$-system with root $r$.
\end{lem}
\begin{proof}
Using $\thry{AC}_\omega(\R)$, fix a $k\in\omega$ such that 
$A:=\{t\in \R: |B_t|=k\}$ is non-meagre. We will have $D\subseteq A$ in the end. 

For each $s\subseteq[\ran B]^{<\omega}$ let $A_s = \{t\in A:  s\subseteq B_t\}$. Let $T$ consist of those $s$ such that $A_s$ is non-meagre, ordered by extension. Notice that $T$ is well-founded: any linearly ordered subset has length at most $k$. Fix a terminal node $r\in T$. Let $n\in\omega,$ $n>1$ be given. We construct $\langle t_i:i\le n\rangle$, eventually letting $D=\{t_i:i\le n\}$.

Fix any $t_0\in A_r$. Notice that if $ s\not\subseteq r$ is finite then $A_{s\cup r}$ is meagre, this will be the key thing to exploit. Assume that $\langle t_i:i< k\rangle$ has been constructed such that $\{B_{t_i}:i<k\}$ forms a $\Delta$-system with root $r$. Notice that $\wp\left(\bigcup\{B_i:i< k\}\right)$ is finite and for each $s\in \wp\left(\bigcup\{B_i:i< k\}\right)$ with $r\subsetneq r\cup s$, $A_{r\cup s}$ is meagre and so $$A_r\bcl\bigcup\{A_{r\cup s}:s\in\wp\left(\bigcup\{B_{t_i}:i<k\}\right)\bcl\wp(r)\}$$ is non-meagre. Fix $t_k$ from this non-meagre set. We show that $\{B_{t_i}:i\le k\}$ remains a $\Delta$-system with root $r$, the only thing to show is that for $i<k$ that $B_{t_i}\cap B_{t_k}=r$ and the choice of $t_k$ makes it obvious that $r\subseteq B_{t_k}\cap B_{t_i}$, if $B_{t_k}\cap B_{t_i}=r\cup s$ for some nonempty $s\subseteq B_{t_i}\bcl r$, then $t_k\in A_{r\cup s}$ which contradicts the choice of $t_k$. 
\end{proof}

\begin{lem}[Miller \cite{miller}]\label{folk}
Suppose a set $G\subseteq\R$ is comeagre in a nonempty open set. $G$ contains a perfect almost disjoint set. 
\end{lem}
\begin{proof} For convenience, $\omega^{<\omega}$ only contains increasing sequences. We will construct a tree $T\subseteq \omega^{<\omega}$ such that whenever $f\in[T]$ then $\ran f\in G$ and whenever $f,g\in[T]$ are distinct then $\ran f\cap \ran g$ is finite. Since $G$ is comeagre in an open neighborhood we may fix a $s\in \omega^{<\omega}$ and a sequence $\langle D_n:n<\omega\rangle$ of open dense-below-$s$ subsets of $\omega^{<\omega}$ such that $D_{n+1}\subseteq D_n$ and $$\bigcap\{\bigcup\{N_s:s\in D_n\}:n<\omega\}\subseteq \hat{G}=\{f_x:x\in G\}$$ where $f_x$ is the increasing enumeration of $x$. We construct $\langle s^i_k:i<2^k, k<\omega\rangle$, $T$ will be the downward closure of $\{s^i_k:k<\omega, i<2^k\}$, we induct on $k$. $\{s^i_k:i<2^k,k<\omega\}$ will be the splitting nodes of $T$ and for any fixed $k$ we will also have $s^i_k\in D_k$ for all $i<2^k$. Let $s_0^0\in D_0$ be any extension of $s$. Once $\{s^i_j:i<2^j,j<k\}$ have been constructed let $n=\max\{\max(\ran s^i_j):i<2^j,j<k\}+1$ we pick for each $s^i_{k-1}$ two extensions such that for any distinct $s^i_{k-1}$ their extensions will have nothing in their range in common above $n$. Let $s_k^0$ be any proper extension of $s_{k-1}^0\ ^\frown n$ in $D_k$ and $s_k^1\in D_k$ any proper extension of $s_{k-1}^0\ ^\frown (\max s_{k}^0+1)$. Once $s^{2i}_k, s^{2i+1}_k\in D_k$ proper extensions of $s^i_{k-1}$ have been defined for all $i<j<2^{k-1}$ for some $j$, let $n' = \max\{\max(\ran s^\ell_k) : \ell<2j\}+1$, $s^{2j}_k\in D_k$ be any proper extension of $s_{k-1}^j\ ^\frown n'$, and $s_k^{2j+1}\in D_k$ any proper extension of $s_{k-1}^j\ ^\frown (\max(s_k^{2j}) + 1).$

That concludes the construction of the tree $T$. It is clear that if $f\in[T]$ then there is a map $\phi_f:\omega\to \omega^{<\omega}$ such that for all $k$, $\phi_f(k)\in\{s^i_k: i<2^k\}$ and $f\restriction |\phi(k)|=\phi(k)$ and so $\{f\}= \bigcap\{N_{\phi_f(k)}:k<\omega\}\subseteq\hat{G}$ and thus for any $f\in [T]$, $\ran f\in G$. Moreover, if $f,g\in [T]$ are distinct then there is some least $k$ such that $\phi_f(k+1)\neq \phi_g(k+1)$ and so $\ran f\cap\ran g\subseteq\ran \phi_f(k)$, which is finite and so $P=\{\ran f:f\in[T]\}\subseteq G$ is a perfect almost disjoint family contained in $G$. 
\end{proof}

To more easily state the lemmas and theorems to disprove the existence of maximally independent families under determinacy-like hypotheses, a definition.

\begin{defn}
Suppose $X$ and $Y$ are sets, $\mathcal{K}$ is an ideal on $X$ and $\sJ_0$ and $\sJ_1$ are ideals on $Y$ with $\sJ_0\subseteq\sJ_1$. Let $\thry{p}\mathcal{K}$ be the positive co-positive pieces, 
$$\thry{p}\mathcal{K}= \{a\subseteq X :a\notin \mathcal{K}\land X\bcl a\notin \mathcal{K}\}.$$

Say $\mathcal{K}\le\langle \sJ_0,\sJ_1\rangle$ if and only if there is a $\Sigma:\thry{p}\mathcal{K}\to\thry{p}\sJ_0$ such that:
\begin{enumerate}
\item[\textup{(1)}]{$\Sigma(\neg a)\triangle \neg\Sigma(a)\in\sJ_0$ for all $a\in\thry{p}\mathcal{K}$ and }
\item[\textup{(2)}]{for all $a,b\in\thry{p}\mathcal{K}$ with $a\cap b\in\mathcal{K}$, $\Sigma(a)\cap\Sigma(b)\in\sJ_1$.}
\end{enumerate}
\end{defn}

\begin{lem}\label{sigmactbl}
Suppose $\kappa$ is a cardinal of countable cofinality, then $\thry{fin}\le\langle \sB(\kappa),[\kappa]^{<\kappa}\rangle$. In addition if $\sJ$ is an ideal on $\kappa$ such that $\sB(\kappa)\subseteq\sJ\subseteq[\kappa]^{<\kappa}$ and $\sF$ is a $\sJ$-independent family on $\kappa$ then $\thry{fin}\le\langle \sJ,\id_{\sJ}(\sF)\rangle$. 
\end{lem}
\begin{proof}
Fix an increasing cofinal $\rho:\omega\to\kappa$. We need to modify $\rho$ slightly to obtain a $\ol{\rho}$. Let $\ol{\rho}(0)=\rho(0)$ and 
$\ol{\rho}(n+1) = \ol{\rho}(n)+\rho(n)$ so that $\rho(n)$ injects into $[\ol\rho(n),\ol\rho(n+1))$ for all $n$, and in particular, if $A$ is a union of infinitely many intervals of the form $[\ol\rho(n),\ol\rho(n+1))$ then $|A|=\kappa$. Define $\Sigma:\R\to\thry{p}[\kappa]^{<\kappa}\subseteq \thry{p}\sJ\subseteq\thry{p}\sB(\kappa)$ by $\Sigma(x) = \bigcup\{[\ol\rho(n),\ol\rho(n+1)):n\in x\}.$ Then $\Sigma(\neg x)=\neg\Sigma(x)$ and if $x,y\in\R$ have $x\cap y$ is finite, then $\Sigma(x)\cap\Sigma(y)$ is bounded and thus $\Sigma(x)\cap\Sigma(y)\in \sB(\kappa)\subseteq\id_\sJ(\sF).$
\end{proof}

\begin{lem}\label{sigmaunctbl}
Suppose there is no injection $\omega_1\to\R$ and $\sJ$ is an $\omega_1$-complete ideal on an ordinal $\kappa$. Then for any infinite $\sJ$-independent family $\sF$, $\thry{fin}\le \langle \sJ,\id_\sJ(\sF)\rangle$. 
\end{lem}
\begin{proof}
Fix a partition $\langle A_n:n\in\omega\rangle$ with each $A_n\in\id_\sJ(\sF)$. Since $\sJ$ is $\omega_1$-complete there must be infinitely many $n$ such that $A_n\in\sJ^+$, otherwise $\kappa\in\id_{\sJ}(\sF)$. Let $\langle n^0_i:i<\omega\rangle$ enumerate those $n$ with $A_n\in\sJ^+$ and $\langle n^1_i:i<k\le\omega\rangle$ enumerate those $n$ such that $A_n\in\sJ$. For all $i<k$ let $\hat{A}_i=A_{n_i^0}\cup A_{n_i^1}$ and if $k<\omega$ for $i\ge k$ let $\hat{A}_{i}=A_{n_i^0}$. So we have a partition $\langle \hat{A}_n:n\in\omega\rangle$ of $\kappa$ such that $\hat{A}_n\in \id_\sJ(\sF)\cap \sJ^+$ for all $n\in\omega$. Then define $\Sigma:\R\to \sJ^+$ by $\Sigma(x) = \bigcup\{\hat{A}_n:n\in x\}.$ It is clear that $\Sigma$ witnesses that $\thry{fin}\le\langle\sJ,\id_\sJ(\sF)\rangle$. 
\end{proof}

\begin{defn}
Suppose $\sJ$ is an ideal on a set $X$ and $a,b\subseteq X$, then we write $a\subseteq_\sJ b$ to mean that $a\bcl b\in\sJ$. 
\end{defn}

\begin{remark} 
Notice that ``$\sR$ is $\sJ$-reaping" is equivalent to ``for all $y\in\sJ^+$ there is an $x\in\sR$ so that either $x\subseteq_\sJ y$ or $x\subseteq_\sJ \neg y.$" Further, notice that $\subseteq_\sJ$ behaves much like $\subseteq $ with respect to finitary operations e.g. if $a\subseteq_\sJ x$ and $b\subseteq_\sJ y$ then $a\cap b\subseteq_\sJ x\cap y$ and so on.
\end{remark}

\begin{thm}\label{main}
Assume $\thry{ZF + DC}$ and all sets of reals have the property of Baire. Suppose $\kappa$ is an ordinal, $\sJ$ is a ideal on $\kappa$, and $\sF$ is an $\sJ$-independent family on $\kappa$ such that $\thry{fin}\le\langle \sJ,\id_\sJ(\sF)\rangle$, then $\sF$ is not maximal.
\end{thm}
\begin{proof} Suppose $\sF$ was maximal, we construct a boolean combination in $\id_\sJ(\sF)$. Since $\sF$ is maximal $\sR:=\sR(\sF)$ is $\sJ$-unsplit and so for each $a\in\sJ^+$, there is some $B\in \sR$ so that either $B\subseteq_\sJ a$ or $B\subseteq_\sJ \kappa\bcl a$. Fix a $\Sigma:\R\to\thry{p}\sJ$ witnessing that $\thry{fin}\le\langle \sJ,\id_\sJ(\sF)\rangle$. Let 
$$G = \{x\in\R : (\exists B\in\sR)(B\subseteq_\sJ \Sigma(x))\}$$ and 
$$H= \{x\in\R:(\exists B\in\sR)(B\subseteq_\sJ \kappa\bcl \Sigma(x))\}.$$
The maximality of $\sF$ implies that $G\cup H=\R$ and since $\Sigma(\neg x)\triangle\neg\Sigma(x)\in\sJ,$ $x\in G\leftrightarrow \neg x\in H$ and so $H$ and $G$ are images of one another under bit-flipping, a homeomorphism of the real line. So, $G$ and $H$ must have the same category, in particular $G$ must be non-meagre. Since $G$ has the property of Baire and is non-meagre we may and do fix a perfect almost disjoint $P\subseteq G,$ appealing to lemma \ref{folk}. Fix injections $\phi:\R\to P$ and $k:[\sF]^{<\omega}\times[\sF]^{<\omega}\to\wp(\kappa)$. Then define a relation $B\subseteq\R\times \wp(\kappa)$ by $\langle t,a\rangle\in B$ if and only if there are $X,Y\subseteq\sF$ finite disjoint such that $k(X,Y)=a$ and $\thry{b}(X,Y)\subseteq_\sJ\Sigma\circ\phi(t).$ Note that $\dom B=\R$. Now we want to choose, for each $t\in\R$, an $a$ such that $\langle t,a\rangle\in B$. Of course there is insufficient choice in $V$, however we only need $V$ to see finitely many of these choices.

Let $g\subseteq\thry{Col}(\omega_1,\wp(\kappa))$ be $V$-generic. Notice that since we are using countable conditions and we have \thry{DC} in $V$ we have $\thry{AC}_\omega(\R)$ in $V[g]$ and $\R^V=\R^{V[g]}$ and so $P$ remains perfect in $V$ and it remains the case that $\dom B=\R$. Then working in $V[g]$ we have well-ordered $\wp(\kappa)^V$. We make our choices, for each $t\in\R$ let $a_t\in\wp(\kappa)^V$ be such that $\langle t,a_t\rangle\in B$ and let $\langle X_t,Y_t\rangle=k^{-1}(a_t)$ so that $\thry{b}(X_t,Y_t)\subseteq_\sJ \Sigma\circ\phi(t).$

Now define a relation $S\subseteq \R\times \wp(\kappa)^V$ by $\dom S=\dom B=\R$ and for all $t\in\R$, $S_t = X_t\cup Y_t.$ Now, we use the $\Delta$-system lemma \ref{delta} applied to $S$ to get a finite $r\subseteq\sF$ and a $D\subseteq\R$ with $|D|>2^{|r|}$ such that $\{S_t:t\in D\}$ forms a $\Delta$-system with root $r$. Each $t\in D$ induces a partition of the root $r$ into two pieces as $\langle r\cap X_t, r\cap Y_t\rangle$. Since $|D|$ is strictly larger than the number of such partitions there must be distinct $t,u\in D$ such that $r\cap X_t=r\cap X_u$ and $r\cap Y_t= r\cap Y_u.$ Now we return to working in $V$, we have $t,u,r,\langle X_t,Y_t\rangle ,\langle X_u,Y_u\rangle\in V.$ Since $\thry{b}(X_t,Y_t)\subseteq_\sJ \Sigma(\phi(t))$ and similarly for $u$ we have that 
$$\thry{b}(X_t,Y_t)\cap\thry{b}(X_u,Y_u)\subseteq_\sJ\Sigma(\phi(t))\cap\Sigma(\phi(u))\in\id_\sJ(\sF)$$ as since $u\neq t$, $\phi(u)\cap \phi(t)\in\thry{fin}$, and so $\Sigma(\phi(u))\cap\Sigma(\phi(t))\in\id_\sJ(\sF)$. However, since $r\cap X_t=r\cap X_u$, $r\cap Y_t= r\cap Y_u$, and $(X_t\cup Y_t)\cap (X_u\cup Y_u)=r$ we have that $(X_t\cup X_u)\cap (Y_t\cup Y_u)=0$ and therefore 
$\thry{b}(X_t\cup X_u, Y_t\cup Y_u)\in\sR(\sF).$ If one then looks closely one notices that 
$$\thry{b}(X_t,Y_t)\cap\thry{b}(X_u,Y_u)= \thry{b}(X_t\cup X_u, Y_t\cup Y_u)$$ and so we have $\sR(\sF)\cap\id_\sJ(\sF)\neq 0$, an absurdity!
\end{proof}

\begin{cor}\label{cort}
Assume $\thry{AD + DC}_\R$ and $\kappa<\Theta$ is a cardinal. Say $\sJ$ is an ideal on $\kappa$ such that either $\sJ$ is $\omega_1$-complete or $\cf(\kappa)=\omega$ and $\sB(\kappa)\subseteq\sJ\subseteq[\kappa]^{<\kappa}$.
There are no infinite maximal $\sJ$-independent families on $\kappa$.
\end{cor}
\begin{proof}
\thry{DC} was used in the previous only to ensure that in a generic extension we have $\thry{AC}_\omega(\R)$ and that reals are not added. Since $\thry{AD + DC}_\R$ holds and $\kappa<\Theta$ there is a surjection $\R\twoheadrightarrow \wp(\kappa)$ thus \thry{DC}$_{\wp(\kappa)}$ holds which is enough for the argument.
\end{proof}

\begin{cor}\label{corm}
Assume $J$ is a set of ordinals and $V=L(\R,J)+\thry{AD + }\thry{DC}_\R$ and $\kappa$ is a cardinal. Say $\sJ$ is an ideal on $\kappa$ such that either $\sJ$ is $\omega_1$-complete or $\cf(\kappa)=\omega$ and $\sB(\kappa)\subseteq\sJ\subseteq[\kappa]^{<\kappa}$.
There is no infinite maximal $\sJ$-independent family on $\kappa$.
\end{cor}

The only possible models satisfying $V=L(\wp(\R))+\thry{AD}^+$ which could have a maximally independent family satisfy $\thry{AD}_\R\thry{ + }\neg\thry{DC}_{\wp(\R)}$. However, if such a model had a maximal independent family then one of its inner models of the form $L(\R)$ would contain an independent family which is sufficiently maximal to not exist.

\begin{thm}\label{lpr}
Suppose $V=L(\wp(\R))\thry{ + AD}^+$ and $\sJ$ is an ideal on a cardinal $\kappa$ such that either $\sJ$ is $\omega_1$-complete or $\cf(\kappa)=\omega$ and $\sB(\kappa)\subseteq\sJ\subseteq[\kappa]^{<\kappa}$.
If $\sF$ is an infinite $\sJ$-independent family on $\kappa$ then $\sF$ is not maximal. 
\end{thm}
\begin{proof}We consider two cases, $\sF$ being countable and $\sF$ being uncountable.  Fix a bijection $\psi:\omega\times\kappa\to\kappa$ with $\psi\in\thry{HOD}.$ We first consider the case that $\sF$ uncountable and $\sJ$-maximal and arrive at a contradiction.  Without loss of generality assume $\thry{AD}_\R\thry{ + }\neg\thry{DC}$ so that $\cf(\Theta)=\omega.$
Since $\thry{AD}^+$ holds the Wadge ranks are well-ordered in type $\Theta$, for $C\subseteq \R$ let $w(C)\in\Theta$ denote the Wadge rank of $C$. We will now construct an inner model $\sM$ in which $\sF$ is ``maximal enough''. For each $\eta<\Theta$ let 
$$\sF_\eta = \{x\in\sF : (\exists C\subseteq\R)(x\in\thry{HOD}_C\text{ and } w(C)\le\eta)\}$$ and $\eta_0\in\Theta$ least with $\sF_{\eta_0}$ infinite, $\eta_0$ exists since $\cf(\Theta)=\omega$ and $\sF$ is uncountable. 

Further, fix a partition $\langle A_n:n<\omega\rangle$ of $\kappa$ into countably many pieces such that, letting $\Sigma:\R\to\thry{p}\sJ$ be defined by $\Sigma(x)=\bigcup\{A_n:n\in x\},$ we have that $\Sigma$ witnesses that $\thry{fin}\le\langle\sJ,\id_\sJ(\sF)\rangle$ as in either lemmas \ref{sigmactbl} or \ref{sigmaunctbl}.
This sequence $\langle A_n:n<\omega\rangle$ is not unique. If we are in the case that $\cf(\kappa)=\omega$ and $\sB(\kappa)\subseteq\sJ\subseteq[\kappa]^{<\kappa}$, then $\langle A_n:n<\omega\rangle$ is defined from a cofinal sequence through $\kappa$ in such a way that $A_n\in\sB(\kappa)\subseteq\id_\sJ(\sF)$ for all $n\in\omega$. In the case that $\sJ$ is $\omega_1$-complete, the sequence $\langle A_n:n<\omega\rangle$ is defined from an injection $\omega\to\sF$. Fix an injection $\vec{x}=\langle x_n:n<\omega\rangle$ of $\omega\to \sF_{\eta_0}$, and let $\langle A_n:n<\omega\rangle$ be defined from $\vec{x}$. That is, $$\{A_n:n<\omega\}=\left\{\bigcap\{x_n:f(n)=0\}\bcl\bigcup\{x_n:f(n)=1 \} : f\in2^\omega\right\}\bcl\{0\}.$$
If $\cf(\kappa)=\omega$ and $\sB(\kappa)\subseteq\sJ\subseteq[\kappa]^{<\kappa}$ then let 
$J=\psi"\{\langle n,\beta\rangle:n\in\omega\land\beta\in A_n\}$
and if else $\sJ$ is $\omega_1$-complete let 
$$J = \psi"\{\langle 2n,\beta\rangle:n\in\omega\land\beta\in A_n\}\cup\psi"\{\langle 2n+1,\beta\rangle:n\in\omega\land\beta\in x_n\}.$$

For each $\eta<\Theta$ let 
$$Z_\eta = \{x\in\R: (\exists B\in\sR(\sF_\eta))(B\subseteq_\sJ\Sigma(x)\lor B\subseteq_{\sJ}\Sigma(\neg x))\}.$$ Since the ideal of meagre sets is $\infty$-complete and $\bigcup\{Z_\eta:\eta<\Theta\}=\R$, there is some least $\eta_1<\Theta$ such that $\eta_0\le\eta_1$ and $Z_{\eta_1}$ is non-meagre. Fix a $C^*\subseteq\R$ with $w(C^*)>\eta_1$ and such that each $\sF_{\eta_1}$ and $\subseteq_\sJ$ are $\thry{OD}_{C}$ for some $C$ with $w(C)\le w(C^*)$.  Let $\sF^*=\sF_{\eta_1}$, $Z^*=Z_{\eta_1}$, and $\mathcal{M} = \thry{HOD}_{\R\cup\{C^*,J\}}.$ Then we have:

\begin{enumerate}[(1)]
\item{$\Sigma$, $\sF^*$, $Z^*$, $\langle A_n:n<\omega\rangle$   
, $\subseteq_\sJ\cap\sM$, and $\sJ\cap\sM$ are all elements of $\sM$}
\item{$\sM\model (\forall n)(A_n\in\id_\sJ(\sF^*))\ \land\ \Sigma \text{ witnesses that }\thry{fin}\le\langle\sJ,\id_\sJ(\sF^*)\rangle$. In the case that $\sJ$ is $\omega_1$-complete this requires that $\{x_i:i<\omega\}\subseteq\sF^*\cap \sM$.  }
\item{$\sM\model\thry{AD + DC + }Z^*$ is non-meagre. Proposition \ref{HODDC} addresses $\sM\model\thry{DC}$.}
\end{enumerate}

We can now essentially run the proof of theorem \ref{main} within $\sM$. Working in $\sM$ we don't necessarily have that $\sF^*$ is maximal, however, if we let $$G^* = \{x\in\R : (\exists B\in \sR(\sF^*))(B\subseteq_\sJ\Sigma(x))\}$$
and 
$$H^* = \{x\in\R : (\exists B\in \sR(\sF^*))(B\subseteq_\sJ\Sigma(\neg x))\}$$
we have that $G^*\cup H^*=Z^*$ is non-meagre and $x\in G^*$ if and only if $\neg x\in H^*$ and so $G^*$ is non-meagre, which is all we need. Then exactly as in the proof of theorem \ref{main} we go into $\sM[g]$ where $g\subseteq\thry{Col}(\omega_1,\wp(\kappa))^{\sM}$ is $\sM$-generic to find $t,u\in\R$ and finite $X_t,Y_t,X_u,Y_u\subseteq\sF^*$ such that $\thry{b}(X_t,Y_t)\subseteq_\sJ \Sigma(t)$, $\thry{b}(X_u,Y_u)\subseteq_\sJ\Sigma(u)$, $t\cap u\in\thry{fin}$, and $(X_t\cup X_u)\cap (Y_t\cup Y_u)=0$. So that in $\sM$ we have a contradiction:
$$\thry{b}(X_t\cup X_u, Y_t\cup Y_u)\in\sR(\sF^*)\cap\id_\sJ(\sF^*).$$ 
The case that $\sF$ is countable is much simpler. By coding the family $\sF$ into a set of ordinals $J$ and finding a $C\subseteq\R$ from which the ideal $\sJ$ is definable one finds that $\sN=\thry{HOD}_{\R\cup\{C,J\}}$ contains as elements each of $\sF, \subseteq_\sJ\cap\sN,$ and $\sJ\cap\sN$ and that $\sN\model\thry{AD + DC + fin}\le\langle\sJ,\id_\sJ(\sF)\rangle$. So, $\sN\model\sF$ is not maximal, and thus $\sF$ is not maximal in $V$ either. 
\end{proof}

Theorem \ref{lpr} completely resolves the problem of the existence of maximal independent
families in $\thry{AD}^+$ models of the form $L(\wp(\R))$ and theorem \ref{main} resolves the same problem in Chang-like models constructed in \thry{AD + DC} models. 

\begin{qstin}
Assume $V=L(\thry{On}^\omega)\thry{ + AD + }\neg\thry{DC}.$ Is there a maximally independent family?
\end{qstin}

\end{document}